\documentclass{article}

\usepackage{proof}
\usepackage{latexsym}
\usepackage{amssymb}

\newcommand{\blem}{\begin{lemma}}
\newcommand{\elem}{\end{lemma}}
\newcommand{\bth}{\begin{theorem}}
\newcommand{\benu}{\begin{enumerate}}
\newcommand{\eenu}{\end{enumerate}}
\newcommand{\bdes}{\begin{description}}
\newcommand{\edes}{\end{description}}
\newcommand{\bdf}{\begin{definition}}
\newcommand{\edf}{\end{definition}}
\newcommand{\bcor}{\begin{cor}}
\newcommand{\ecor}{\end{cor}}

\newcommand{\beqn}{\begin{equation}}
\newcommand{\eeqn}{\end{equation}}

\newtheorem{theorem}{Theorem}
\newtheorem{definition}{Definition}

\newtheorem{lemma}{Lemma}
\newtheorem{cor}{Corollary}

\newcommand{\alp}{\alpha}
\newcommand{\veps}{\varepsilon}
\newcommand{\del}{\delta}
\newcommand{\Del}{\Delta}
\newcommand{\ome}{\omega}
\newcommand{\Ome}{\Omega}
\newcommand{\bet}{\beta}
\newcommand{\gam}{\gamma}
\newcommand{\Gam}{\Gamma}
\newcommand{\kap}{\kappa}
\newcommand{\sig}{\sigma}
\newcommand{\Sig}{\Sigma}
\newcommand{\lam}{\lambda}

\newcommand{\vphi}{\varphi}
\newcommand{\fal}{\forall}
\newcommand{\exi}{\exists}
\newcommand{\rarw }{\rightarrow}
\newcommand{\Rarw }{\Rightarrow}

\newcommand{\Lrarw}{\Leftrightarrow}

\newcommand{\calr}{{\cal R}}
\newcommand{\calL}{{\cal L}}
\newcommand{\msten}{\mbox{\hspace{10mm}}}

\title{Introducing the Hardline in Proof Theory
\thanks{This had been submitted to the Bulletin of Symbolic Logic as a communication in May 1996. According to the referee's report received Sep.1996, even for an expert the paper is too sketchy and only a very small part is accessible to a wide audience. The editor decided not to accept this for publication in the Bulletin. The hardliner withdrew this from publication.}}
\author{Toshiyasu Arai\\
Faculty of Integrated Arts and Sciences\\
Hiroshima University
\thanks{current address:
Graduate School of Science,
Chiba University,
1-33, Yayoi-cho, Inage-ku,
Chiba, 263-8522, JAPAN}
}
\date{}
\begin{document}
\maketitle

G. Gentzen \cite{G3} published his new version of consistency proof for first order number theory in 1938. He had already had two consistency proofs \cite{G1} and \cite{G2}. The first used {\em constructive but rather abstract notion of functionals}. In the second he had first introduced {\em transfinite ordinals} in proof theory. Although he formulated the result as a consistency proof, his interest seems to involve a taking off from Hilbert's program. As to this turning G. Kreisel \cite{K} p. 262 wrote: 

\smallskip
\begin{quotation}
$\ldots$, by introducing a {\em quantitative ordinal measure} he (=Gentzen) forces us to pay attention to {\em combinatorial complexity}\footnote{This emphasis is made by the author.} and thereby makes it at least more difficult for us to slip into an abstract reading.
\end{quotation}
\smallskip

It seems that the purpose of the third "Neue Fassung" is to make a lucid exposure of this combinatorial complexity which Gentzen discovered in {\em finite proof figures} of number theory.

G. Takeuti followed this idea and developed a proof theory of systems of second order arithmetic including $\Pi^1_1$-Comprehension Axiom, $\Pi^1_1-CA$.\\
We follow in the wake of Gentzen and Takeuti. Proof theory \`a la Gentzen prodeeds as follows;
\bdes
\item[(G1)] Let $P$ be a proof whose endsequent has a restricted form. Define a reduction procedure $r$ which rewrites such a proof $P$ to yield another proof $r(P)$ provided that $P$ has not yet reduced to a certain canonical form.
\item[(G2)] From the structure of the proof $P$, we abstract a structure related to this procedure $r$ and throw irrelevant residue away. Thus we get a finite figure $o(P)$.

We call the figure $o(P)$ the {\em ordinal diagram} ( abbr. by o.d.'s) after G. Takeuti \cite{T}. Let ${\cal O}$ denote the set of o.d.'s.
\item[(G3)] Define a relation $<$ on ${\cal O}$ so that $o(r(P))<o(P)$.
\item[(G4)] Show the relation $<$ on ${\cal O}$ is well founded. \\
Usually $<$ is a linear ordering and hence $({\cal O},<)$ is a notation system for ordinals.
\edes

This description is not acute. In fact (G1)-(G4) interact each other. For example (G1) is influenced by (G3) and this by (G4).

In this paper we expound some basic ideas of proof theory for theories of ordinals $\sig$ such that there are many $\sig$-stable ordinals below $\sig$. From this we get the proof theoretic ordinals of subsystems of second order arithmetic, e.g., $\Sig^1_3-DC+BI$. The deatils will be reported in \cite{A2}, \cite{A3},\cite{A4}.

In \S 2 we expound proof theory for $\Pi_3$-reflecting and $\Pi_4$-reflecting ordinals in some detail. In \S 3 theories for ordinals $\sig$ having many $\sig$-stable ordinals below are analysed.

For more on the aims and another approach to proof theory of strong theories, see M. Rathjen \cite{R1} and \cite{R2}.

\section{$\Pi_2^\Ome$-ordinal of a theory}
    G. J\"ager \cite{J} has shifted an object of proof-theoretic study to set theories from second order arithmetic.

\begin{definition}($\Pi_2^\Ome$-ordinal of a theory)  
{\rm Let} $T$ {\rm be a recursive theory of sets such that}
 $KP\omega\subseteq T\subseteq ZF+V=L$, {\rm where} $KP\omega$ {\rm denotes Kripke-Platek set theory with the Axiom of Infinity.} {\rm For a sentence} $A$ {\rm let}  $A^{L_\alpha}$ {\rm denote the result of replacing unbounded quantifiers} $Qx \, (Q\in\{ \forall, \exists\})$ {\rm in} $A$ {\rm by} $Qx\in L_\alpha$. {\rm Here for an ordinal} $\alpha\in Ord$ $L_\alpha$  {\rm denotes an intial segment of G\"odel's constructible sets.} {\rm Let} $\Omega$ {\rm denote the (individual constant corresponding to the) ordinal} $\omega^{CK}_1$. {\rm If} $T\not\vdash\exists \omega^{CK}_1$, {\rm e.g.,} $T=KP\omega$ , {\rm then} $A^{L_\Ome }=_{df}A $. {\rm Define the}
$\Pi_2^\Ome$-ordinal $\mid T\mid$ of  $T$ {\rm by}
\[\mid T\mid=_{df}\inf\{\alpha\leq\omega^{CK}_1 :\forall\Pi_2\, sentence\, A(T\vdash A^{L_\Ome} \: \Rightarrow \: L_\alpha\models A)\}<\omega^{CK}_1\]
\end{definition}
Here note that $\mid T\mid<\omega^{CK}_1$ since we have
\[\forall\Pi_2\, sentence\, A(T\vdash A^{L_\Ome} \: \Rightarrow \: L_\Ome\models A)\]
and $\Ome=\ome^{CK}_1$ is recursively regular, i.e., $\Pi_2$-reflecting.

G. J\"ager \cite{J} shows that $\mid KP\omega\mid$=Howard ordinal and G. J\"ager and W. Pohlers \cite{J-P} gives the ordinal $\mid KPi\mid$, where $KPi$ denotes a set theory for recursively inaccessible universes. Also see Buchholz and Sch\"utte \cite{B-S} and Sch\"utte \cite{S} for related results. These include and imply proof-theoretic ordinals of subsystems of second order arithmetic corresponding to set theories. 

We will develop proof theory for theories of ordinals:
Let ${\cal L}_0$ denote the first order language whose constants are;
$=$(equal), $<$(less than), $0$(zero), $1$(one), $+$(plus), $\cdot$(times), $j $(pairing),$( )_0 ,( )_1$(projections,i.e., inverses to $j$).
\\
For each $\Delta_0$ (=bounded) formula  ${\cal A}(X,a,b)$ (a binary predicate $X$ ) we introduce a binary constant $R^{\cal A}$ such that
\[b\in R^{\cal A} \: \Leftrightarrow_{df} \:  R^{\cal A} (a,b)\:  \Leftrightarrow  \: {\cal A}(R^{\cal A}_{<a} ,a,b) \]
  with $R^{\cal A}_{<a} =\sum_{x<a}R^{\cal A}_x = \{(x,y):x<a \, \& \, y\in R^{\cal A}_x\}$.
\\
Let $F:Ord\rightarrow L$ denote (a variant of) the G\"odel's enumeration of constructible sets.
Then one can define relations $\varepsilon$ and $\equiv$ on $Ord$ such that 
\[\alpha\varepsilon\beta \:  \Leftrightarrow \: F'\alpha\in F'\beta \: ; \: \alpha\equiv\beta \: \Leftrightarrow \: F'\alpha=F'\beta.\] 
and these relations $\varepsilon$ and $\equiv$ are definable by $\Del_0$ fomulae in the language $L_0\cup\{R^{\cal A}\}$.
\medskip

Thus, in principle, one can define a theory $T^{Ord}$ of ordinals for each set theory $T$  by interpreting $\in$ and $=$ as $\varepsilon$ and $\equiv$, resp. In place of $T ^{Ord}$ we consider a theory $T_n$ of $\Pi_n$-reflecting ordinals.

\begin{definition} (Aczel and Richter \cite{A-R})  
{\rm Let} $X\subseteq Ord$ {\rm denote a class of ordinals and} $\Phi$ {\rm a set of formulae in the language of set theory (or the language of theories of ordinals). Put} $X\!\mid\!\alpha=_{df}\{\beta\in X :\beta<\alpha \}$. {\rm We say that an ordinal} $\alpha\in Ord$ {\rm is} $\Phi$-reflecting on $X$  {\rm if} 
\[\forall A \in\Phi \mbox{{\rm with parameters from }} L_\alpha [L_\alpha\models A \: \Rightarrow \: 
\exists\beta\in X\!\mid\!\alpha (L_\beta\models A)]\]
 {\rm If a parameter} $\gamma<\alpha$ {\rm occurs in} $A ${\rm , then it should be understood that} $\gamma<\beta${\rm .}
\\
$\alp$ {\rm is} $\Phi$-reflecting {\rm if} $\alpha$ {\rm is} $\Phi${\rm -reflecting on the class of ordinals} $Ord$.
\end{definition}

\section{$\Pi_3$ and $\Pi_4$ reflection}

Let us explain how to design a notation system $O(\Pi_3)$ of ordinals (its elements are called {\em ordinal diagrams} abbr. by o.d.'s) and show
\[\forall\Pi_2\: A(T_3\vdash A \: \Rightarrow \: \exists\alpha\in O(\pi_3)\!\mid\!\Omega(=\{\alpha\in O(\pi_3):\alpha<\Omega\}) \: s.t.\: A^\alpha ).\]
@@@
    $T_3$ is formulated in Tait's logic calculus, i.e., one-sided sequent calculus and $\Gamma,\Delta\ldots$ denote a {\em sequent}, i.e., a finite set of formulae. $T_3$ has the inference rule $(\Pi_3-rfl)$:

    \[\infer[(\Pi_3-rfl)] {\Gamma}{\Gamma,A & \neg A^b,\Gamma}\]                          
where $A\equiv\forall x\exists y\fal zB$ with a bounded formula $B$ and the eigenvariable $b$.
\\
So $(\Pi_3-rfl)$ says   $A\rarw \exists bA^b$.\footnote{For simplicity we suppress the parameter. Correctly $\forall u(A(u) \, \rarw \, \exists z(u<z \, \& \,A^z(u)))$.}

To deal with the rule $(\Pi_3-rfl)$ we introduce a new rule:

\[\infer[(cp)]{\Gam,A^{\alp_0}}{\Gam,A}\]
where $A $ is a $\Pi_3$-sentence as above.
\\
\smallskip
We need to compute an o.d. $\alp_0<\pi$ in order to replace the $(\Pi_3-rfl)$ by a $(cut)$:

\[\infer[(cut)]{\Gam}
{\infer[(cp)]
 {\Gam,A^{\alp_0}}{\Gam,A}
&
\deduce{\neg A^{\alp_0},\Gam}{[b:=\alp_0]}
}\]

Firstly we throw $0$ and $\pi$ into $O(\Pi_3)$. The o.d. $\pi$ correspods to the first $\Pi_3$-rfl ordinal. Let $O(\Pi_3)$ be closed under $+$ and the Veblen function $\varphi$. The Veblen function $\varphi$ is needed for treating the constant $R^{\cal A}$. Let $\calr$ denote the set of o.d.'s corresponding to recursively regular ordinals.

We have learnt the following fact from the proof theory for the universes with many recursively regular ordinals:
In general, if $\sig$ is recursively regular, then we have to introduce a collapsing 
$(\sig,\alp)\mapsto d_\sig\alp$. \\
For example, it suffices to have two steps collapsings for recursively Mahlo ordinals: \\
$(\mu,\alp)\mapsto d_\mu\alp=\sig$ and $(\sig,\bet)\mapsto d_\sig\bet $ with the first recursively Mahlo ordinal $\mu$.
\\
    The relation  $\alp<\bet$ is defined so as to hold:
\bdes
\item [($<1$)] $d_\sig\alp<\sig$
\item [($<2$)]$K_\sig\alp<d_\sig\alp$
\item [($<3$)]$K_\sig\alp\leq\alp$
\item[($<4$)] $\alp<\sig \, \&\, K_\sig\alp<d_\sig\bet \: \Rarw \: \alp<d_\sig\bet$
\edes
where $K_\sig\alp$ denotes the finite set of subdiagrams $\bet$ of $\alp$ such that, in the construction of $\alp$, $\bet$ is a last collapse of $\sig$, i.e., 
\[\exi \{\sig_i\}_{i\leq n}\fal i<n(\sig=\sig_0 \, \& \, \sig_{i+1}=d_{\sig_i} \, \& \, \sig_n=\bet)\]

The first candidate to $\alp_0$ is $d_\pi\alp$ with $\alp=o(\Gam,A)$, where $o(\Gam)$ denotes the o.d. assigned to the sequent $\Gam$. But this does not work.
Consider a proof with nested rules $(\Pi_3-rfl)\,J,J_1$:

\[\infer[J]{}
{\infer*{A}
 {\infer[ J_1]{}{A_1 & \neg A_1^{b_1}}
 }
&
\neg A^b
}\]

$A_1\equiv\fal x_1\exi y_1\fal z_1B_1, \: A\equiv\fal x\exi y\fal zB$.
\\
First replace the lower $(\Pi_3-rfl)\, J$ by a $(cp) \, K_0$ followed by a $(cut) \, I$:

\[\infer[I]{}
{\infer[K_0]{A^{\alp_0}}
 {\infer*{A}
  {\infer[J_1]{}{A_1 & \neg A_1^{b_1}}
   }
  }
&
\deduce{\neg A^{\alp_0}}{[b:=\alp_0]}
}
\msten  Fig.1\]
with an o.d. $\alp_0<\pi$, e.g., $\alp_0=d_\pi\alp$.
\\
Then do the same thing to the above $(\Pi_3-rfl)\, J_1$:

\[\infer{}
{\infer[(cp)]{A_1^{\alp_1}}{A_1}
&
\deduce{\neg A_1^{\alp_1}}{[b_1:=\alp_1]}
}\]

We are forced to have $\alp_1<\alp_0$ since $\alp_1$ may be substituted for $y$ in $\exi y\fal zB$, i.e., $\exi y<\alp_0\fal z<\alp_0B$.  But the innermost unbdd universal quantifier $\fal z$ in $A$ causes troubles since any o.d. $\bet<\alp_0$ may be substituted for $z$, e.g., $\bet\geq\alp_1$, and this destroies the case

\[\infer[(cp)]{\exi y_1<\alp_1\fal z_1<\alp_1B_1}
{\infer[(\exi)]{\exi y_1\fal z_1B_1(y)}{\fal z_1B_1(\bet)}
}\]

We cannot anticipate that what o.d. $\bet$ is substituted for $z$ except $\bet<\alp_0$ and $\bet$ comes from the right upper part of the $(cut)\, I$.\\
    How to get rid of this difficulty? Our answer is to {\em iterate collapsings}: Put $\alp_0=d_\pi\alp$ and a $(cp), K_1$ resolving the $(\Pi_3-rfl)\, J_1$ situates below $I$ :

\[\infer[(cut)]{}
{
 \infer[K_1]{A_1^{\alp_1}}
  {
   \infer[I]{A_1^{\alp_0}}
    {
     \infer[K_0']{A^{\alp_0},A_1^{\alp_0}}
      {
       \infer*{A,A_1}{A_1}
       }
    &
    \infer* {\neg A^{\alp_0}}{\bet}
    }
   }
 &
  \deduce{\neg A_1^{\alp_1}}{[b_1:=\alp_1]}
} 
\msten Fig.2\]
with $\alp_1=d_{\alp_0}\bet_1, \: \bet_1=o(A_1^{\alp_0})$

Then  $\bet<\alp_0 \, \&\, K_{\alp_0}\bet<\alp_1 \, \Rarw \, \bet<\alp_1$ is seen from ($<4$). $K_{\alp_0}\bet<\alp_1$ is satisfied since $K_{\alp_0}\bet\subseteq K_{\alp_0}\bet_1<d_{\alp_0}\bet=\alp_1$ by ($<2$).\footnote{$d_\pi\alp\in\calr$ since, in general, the closure ordinal $\bet$ is recursively regular with $A^\pi \Rarw \exi\bet<\pi A^\bet$ for a $\Pi_3 \: A$, cf. \cite{A-R}.}

    In this way we reduce proof figures. The problem is that we have an infinite iteration of collapsings in $O(\pi_3)$: $\pi\mapsto d_\pi\mapsto d_{d_\pi}\mapsto\cdots$. Thus we have readily an infinite decreasing sequence by the requirement ($<1$) $d_\sig\alp<\sig$. We have to kill this infinite sequence. Let us examine what changes when we pass from $\alp_0$ to $\alp_1$ .

    Observe that the upper part of the $(cp)\, K_0$ in $Fig.1$ becomes simpler in $Fig.2$, i.e., the $(cp)\, K_0'$ . This reflects to o.d.'s so that 
 $o(A,A _1) =o(K_0' ) < o(K_0) =o(A)$.  Therefore when we iterate collapsings, i.e., build a tower of rules $(cp)$ growing downwards, the upper part of the topmost $(cp)$ becomes simpler, i.e., o.d. decreases because of resolving $(\Pi_3-rfl)$. Hence when we introduce an o.d. $d_\sig\alp$ from $(\sig,\alp)$ we attach the o.d. $\mu$ corresponding to this upper part: $(\sig,\alp,\mu)\mapsto d_\sig^\mu\alp$. We call the o.d. $\mu$ the {\em q-part} of the o.d. $d_\sig^\mu\alp$ and denote $\mu=st(d_\sig^\mu\alp)$.\footnote{$st$ stands for {\em Stufe} or {\em stage}.} And require that:
\beqn\label{eqn:1}
st(d_\sig^\mu\alp)<st(\sig)  \mbox{ if } \sig\neq\pi, \mbox{ i.e., } \mu<\nu \mbox{ for } d^\mu_{d^\nu_\pi\bet}\alp
\eeqn

Then it may be the case that any infinite collapsing processes are killed by this proviso (\ref{eqn:1}).

    Nonetheless this is not the end of the story. First $\mu=st(d_\sig^\mu\alp)\geq\pi$ in general and so a well ordering proof may be difficult. Further, on the side of proof figures, the proviso (\ref{eqn:1}) means that we have to pinpoint, for each $(cp)$, the unique succession of rules $(cp)$, called the {\em chain}, which describes how to introduce the $(cp)$: For each 
    \[\infer[(cp)^\sig_{d_\sig^\mu\alp}]{A^{d_\sig^\mu\alp}}{A^\sig},\]
     pinpoint the unique chain 
\[\infer[(cp)^\pi_{\sig_1}]{A^{\sig_1}}{A^\pi}, \;    \infer[(cp)^{\sig_1}_{\sig_2}]{A^{\sig_2}}{A^{\sig_1}},\ldots, \infer[(cp)^\sig]{A^{d_\sig^\mu\alp}}{A^\sig}\]
 such that $\sig_1=d_\pi^{\mu_1}\alp_1, \sig_2= d_{\sig_1} ^{\mu_2}\alp_2,\ldots$
  
These $(cp)$'s are connected or related each other by collapsing. And furthermore it must be the case  $o(\mbox{the upper part of the topmost } (cp)^\pi )\leq\mu$, and this topmost $(cp)^\pi$  must be determined uniquely from each $(cp)^\sig$. For otherwise suppose there are two chains for a $(cp)$:

\[\infer[(cp)^\sig_{d_\sig^\mu\alp}]{}
{
 \infer* [chains]{A^\sig}
  {\infer[top \: J_0]{A^{\sig_1}}{\deduce{A^\pi}{(\Pi_3-rfl)\:I_0}}
  &
   \infer[top \: J_1]{ A^{\sig_1}}{\deduce{A^\pi}{(\Pi_3-rfl)\:I_1}}
  }
}
\]
$J_0, J_1$  are topmost ones of chains. Even if we have  $o(J_0),o(J_1)\leq\mu$, there may be $(\Pi_3-rfl)$'s $I_0$ and $I_1$ above $J_0$ and $J_1$, resp.
Here we cannot anticipate which one of
$o(J_0)$ and $o(J_1)$ is bigger. So the proviso (\ref{eqn:1}) breakes down.\\
    To retain the uniqueness of the chain, i.e., not to branch or split a chain, we have to be careful in resolving rules with two uppersequents. \\
     Let us examine more closely the situation since this is instructive for $\Pi_4-rfl$. Our guiding principles are:
\bdes
\item[(ch1)] For any $\infer[(cp)^\sig_\tau]{A^\tau}{A^\sig}$ with $\tau=d_\sig^\mu\alp$, if an o.d. $\bet$ is substituted for an existential quantifier $\exi y<\sig$ in $A^\sig$, i.e., $\bet$ is a realization for $\exi y<\sig$, then $\bet<\tau$, and
\item [(ch2)] Resolving rules such as $(cut)$ must not branch a chain.
\edes
\smallskip
1)  First resolve a $(\Pi_3-rfl)$:

\[\infer[J_0\,(cut)]{}
 {
  \infer[(cp)^\pi_\sig]{A^\sig}{A} 
 &
 \neg A^\sig
 }
\msten  Fig.3\]
with $A\equiv \fal x_1\exi x_2\fal x_3A_3, \, \sig=d_\pi^\mu\alp$.\\
Then resolve the $(cut) \, J_0$:

\[\infer[J_1\, (cut)]{}
 {
  \infer{\neg A_1^\sig}
   {
    \infer{A^\sig}{A}
   &
    \neg A^\sig,\neg A_1^\sig
   }
 &
  \infer[I_0]{A_1^\sig}{A_1}
 }
\msten Fig.4\]
with a $\Sig_2 \: A_1$.
\\
2)  Second resolve a $(\Pi_3-rfl)$ above the $(cp)\, I_0$ and a $(cut)$ as in 1):
 
\[\infer{}
 {
  \infer{\neg B_1^\tau}
   {
     \infer{B^\tau}
      {
        \infer[(cut)]{B^\sig}
         { 
           \neg A_1^\sig        
          &
            \infer[\tilde{I}_0\, (cp)^\pi_\sig]{A_1^\sig,B^\sig}
                                   {\deduce[P_5]{A_1,B}{}}
          }
       }
    &
      \infer{\neg B^\tau,\neg B_1^\tau}
        {
           \neg A_1^\sig
         &
           \infer{A_1^\sig,\neg B^\tau,\neg B_1^\tau}
                 {A_1,\neg B^\tau,\neg B_1^\tau}
         }
    }
 &
  \infer[(cp)^\sig_\tau]{B_1^\tau}
    {
      \infer[J_1]{B_1^\sig}
        {
         \neg A_1^\sig
         &
          \infer[(cp)^\pi_\sig]{A_1^\sig,B_1^\sig}{A_1,B_1}
         }
     }
 }
\]
\[\mbox{\hspace{100mm}} Fig. 5\]

with $\tau=d_\sig^\nu\bet$, a $\Sig_2 \: B_1\equiv\exi y_2\fal y_3B_3$.
 \\          
After that resolve the $(cut)\, J_1$:
\samepage
{
\[\infer{}
 {
  \infer{\neg B_1^\tau}
   {
     \infer{B^\tau}
      {
        \infer{B^\sig}
         { 
           \neg A_1^\sig        
          &
            \infer[\tilde{I}_0]{A_1^\sig,B^\sig}
                                   {\deduce[P_6]{A_1,B}{}}
          }
        }
    &
      \neg B^\tau,\neg B_1^\tau
    }
 &
  \infer{B_1^\tau}
    {
      \infer{B_1^\sig}
        {
           \infer{B_1^\sig,A_2^\sig}
            {
               \neg A_1^\sig
             &
               \infer{A_1^\sig,B_1^\sig,A_2^\sig}{A_1,B_1,A_2}
             }
         &
            \infer[J_0']{\neg A_2^\sig}
             {
                \infer{A^\sig}{A}
              &
                \neg A^\sig,\neg A_2^\sig
              }
         }
     }
 }
\]
\[\mbox{\hspace{100mm}} Fig. 6\]
}
Then resolve the $(cut) \, J_0'$ :

\[\infer[K]{}
 {
  \infer{\neg B_1^\tau}
   {
     \infer{B^\tau}
      {
        \infer{B^\sig}
         { 
           \neg A_1^\sig        
          &
            \infer[\tilde{I}_0]{A_1^\sig,B^\sig}
                                   {\deduce[P_7]{A_1,B}{}}
          }
       }
    &
      \neg B^\tau,\neg B_1^\tau
    }
 &
  \infer[I_1]{B_1^\tau}
    {
      \infer[J_2\, (cut)]{B_1^\sig}
        {
           B_1^\sig,A_2^\sig
         &
            \infer{\neg A_2^\sig}
             {
                \neg A_2^\sig,\neg \tilde{A}_1^\sig
              &
                \infer[(cp)^\pi_\sig\, I_0']{\tilde{A}_1^\sig}
                  {
                    \infer*{\tilde{A}_1}{(\Pi_3-rfl) \, H}
                  }
              }
         }
     }
 }
\; Fig.7\]
3) Thirdly assume that we resolve a $(\Pi_3-rfl) \, H$ above the $(cp)^\pi_\sig \, I_0'$. We introduce a new $(cp)^\sig_\rho \, I_1'$ with $\rho=d_\sig^\eta\gam$ immediately above the $(cut) \,J_2$. Then the new $(cp)^\sig_\rho \, I_1'$ is introduced after the $(cp)^\sig_\tau\,  I_1$  and so $\rho= d_\sig^\eta\gam<\tau$. Hence a new $(cut) \, K'$ is introduced below the $(cut) \, K$:

\[\infer[K']{}
{
  \infer[K]{D^\rho}
   {
     \neg B_1^\tau
   &
     \infer[(cp)^\sig_\tau]{B_1^\tau,D^\rho}
    {
      \infer[J_2]{B_1^\sig,D^\rho}
        {
           B_1^\sig,A_2^\sig
         &
            \infer[I_1 \, (cp)^\sig_\rho]{\neg A_2^\sig,D^\rho}
             {
               \infer{\neg A_2^\sig,D^\sig}
                 {
                   \neg A_2^\sig,\neg \tilde{A}_1^\sig
                 &
                   \infer[(cp)^\pi_\sig]{\tilde{A}_1^\sig,D^\sig}
                         {\tilde{A}_1,D}
                  }
               }
           }
       }
     }
 &
 \neg D^\rho
 }
\msten Fig.8\]
with $D\equiv\fal z_1\exi z_2\fal z_3D_3$.
\\
The principle (ch1) will be retained for the $(cp)^\sig_\rho\, I_1'$ since $\neg A_2$ is a $\Sig_1$ sentence. The principle (ch2) is retained when the $(cut) \, J_2$ is resolved: $A_3$ is a bounded formula and so $A_3^\sig\equiv A_3$. $\neg A_3$ exists above the $(cp)^\sig_\rho\,I_1'$. Therefore the {\em grade} $gr(A_3)$ of the formula $A_3$ which is determined from o.d.'s$<\rho$ occurring in $A_3$ is $gr(A_3)<\rho$.  Thus the new $(cut)$ with the cut formula $A_3$ is introduced below the $(cut) \, K'$.
\\
\smallskip
\\
4)  Next consider the $\Pi_4-rfl$. Assume that $A_3\equiv\exi x_4A_4$ in the above figures. Then one cannot resolve the $(\Pi_4-rfl) \, H$ above the $(cp)^\pi_\sig \, I_0'$ by introducing a $(cp)^\sig_\rho$ with $\rho<\tau$ and a $(cut)$ of the cut formula $D^\rho$. This is seen as in $\Pi_3-rfl$, i.e., because $\neg A_2$ is a $\Sig_2$ sentence. Therefore the chain for $H$ have to connect or merge with the chain $I_0-I_1$ for $B$:
 
 \[\infer[(cp)^\tau_\rho\, I_2]{D^\rho}
   {
    \infer{D^\tau}
     {
      \infer{\neg B_1^\tau}
      {
        \infer[I_1']{B^\tau}
         {
          \infer{B^\sig}
           { 
            \neg A_1^\sig        
           &
            \infer[\tilde{I}_0]{A_1^\sig,B^\sig}
                                   {\deduce[P_9]{A_1,B}{}}
            }
          }
       &
         \neg B^\tau,\neg B_1^\tau
       }
      &
       \infer[I_1]{B_1^\tau,D^\tau}
       {
        \infer{B_1^\sig,D^\sig}
        {
           B_1^\sig,A_2^\sig
         &
           \infer{\neg A_2^\sig,D^\sig}
             {
                \infer{\neg A_2^\sig,\neg \tilde{A}_1^\sig}
                 {
                   \infer[I_0"]{A^\sig}{A}
                  &
                   \neg A^\sig, \neg A_2^\sig,\neg \tilde{A}_1^\sig
                  }
              &
                \infer[I_0']{\tilde{A}_1^\sig,D^\sig}
                                         {\tilde{A}_1,D}
              }
          }
        }
      }
    }
\]
\[\mbox{\hspace{100mm}} Fig. 9\]
with $\rho=d_\tau^\eta\gam$ and a $(cut)$ with the cut formula $D^\rho$ follows this figure as in $Fig.8$.

Then the principle (ch1) for the new $(cp)^\tau_\rho \, I_2$ will be retained similarly for $\Pi_3-rfl$. The problem is that the proviso (\ref{eqn:1}) for $O(\Pi_3)$ may break down; it may be the case $\nu=st(\tau)\leq st(\rho)=\eta$ since we cannot expect the upper part of $(cp)^\pi_\sig\, I_0'$ is simpler than the one of $(cp)^\pi_\sig\, I_0$. 
\\
In other words a new succession $I_0'-I_1-I_2$ of collapsings starts. If this chain $I_0'-I_1-I_2$ would grow downwards as in $\Pi_3-rfl$, i.e., in a chain 
$I_0'-I_1-I_2-\cdots-I_n$, $I_n$ would come only from the upper part of $I_0'$, then the proviso (\ref{eqn:1}) would suffice to kill this process. But the whole process may be iterated : in $Fig.9$ another succession $I_0"-I_1-I_2-I_3$ may arise by resolving the $(cut) \, J_0'$ with a $\Pi_4$ cut formula.
\\
    Nevertheless still we can find a reducing part, that is, the upper part of the $(cp)^\sig_\tau \, I_1$: the upper part of the $(cp)^\sig_\tau \, I_1$ becomes simpler in the step $I_2-I_3$. Therefore in $O(\Pi_4)$ the $q$-part of an o.d. consists of two factors: \[(\tau,\alp,\eta,\pi,\nu,\sig)\mapsto d_\tau^{\eta\pi\nu\sig}\alp=\rho.\]
We set:
\[rg_4(\rho)=\pi,st_4(\rho)=\eta,rg_3(\rho)=\sig,st_3(\rho)=\nu .\]
$\nu=st_3(\rho)$ corresponds to the upper part of a $(cp)^\sig$ while $\sig=rg_3(\rho)$ indicates that the merging point for a chain ending with a $(cp)^\tau_\rho$ is a rule $(cp)^\sig$. \\
    Now the provisos for $O(\Pi_4)$ run as follows:
\beqn\label{eqn:2}
   \mbox{For }  \rho=d_\sig^{\mu\pi}\alp, \; \mu=st_4(\rho)<st_4(\sig)
\eeqn
This corresponds to the case when a $(cp)^\sig_\rho$ is introduced as a resolvent of a $(\Pi_4-rfl)$ above the top of the chain whose bottom is a $(cp)_\sig$.
\beqn\label{eqn:3}
   \mbox{For } \rho=d_\sig^{\eta\pi\nu\sig}\alp, \; \nu=st_3(\rho)<st_3(\kap) 
\eeqn
,where $\kap$ denotes the longest o.d. $\kap\geq\tau$ such that $rg_3(\kap)=\sig$ and $\kap$ is a suffix of a $d$ in $\rho$, e.g., $\kap=\tau$ or $\tau=d_\kap^{-}\bet$, etc.
\\
This corresponds to the case when a $(cp)^\tau_\rho$ is introduced with a merging point $(cp)^\sig$ and previously a $(cp)_\kap$ was introduced with the same merging point $(cp)^\sig$.
\\
\smallskip
{\bf Remark}.  In fact we have a stronger relation $st_3(\rho)\ll_{\sig^+}st_3(\kap)$ rather than mere $st_3(\rho)< st_3(\kap)$,  and this is needed for a well ordering proof.
\\
\smallskip

    Let us try to prove that there is no infinite succession $\pi=\sig_0,\sig_1,\ldots$ of collapsing with $\sig_{n+1}=d_{\sig_n}$. Assume such an infinite sequence exists. It suffices to show, then, there would exist an infinite subsequence $\{\sig_{n_i}\}_{i\in\ome}$ such that 
\[\fal i\in\ome[st_4(\sig_{n_{i+1}})<st_4(\sig_{n_i})]\]

Such a subsequence $\{\sig_{n_i}\}$ ammounts to a subseries $\{I_{n_i}\}$ of the infinite chain $\{I_n\}$ such that each $I_{n_i}$ is introduced as a resolvent of a $(\Pi_4-rfl)$ above $I_{n_0}$.\\
Consider the case when
\[\exi\tau[\#\{n\in\ome:rg_3(\sig_n)=\tau\}=\aleph_0], \, i.e., \: \exi\{\sig_{n_i}\}\fal i\in\ome[ rg_3(\sig_{n_i})=\tau]\]
Then by the proviso (\ref{eqn:3}) we would have
\[\fal i\in\ome[st_3(\sig_{n_{i+1}})<st_3(\sig_{n_i})]\]
We can expect this is not the case. And what else? There may be the case 
\[\fal\tau[\#\{n\in\ome:rg_3(\sig_n)=\tau\}<\aleph_0]\]
This means that the new merging points go downwards unlimitedly. For example in $Fig.9$ a new succession with a merging point $(cp)^\tau_\rho\, I_2$ arises by resolving a $(cut)$ below the $(cp)^\sig_\tau \, I_1'$, i.e., $\tilde{I}_0-I_1'-I_2-I_3 \, (cp)^\rho_\kap$ with a $\kap=d_\rho^{\lam\pi\xi\tau}\del$. But in this case we have
  \[\lam=st_4(\kap)<st_4(\tau)=\nu\]
$st_4(\kap)$ corresponds to the upper part $P_5$ of a $(cp)^\pi_\sig \, \tilde{I}_0$ in $Fig.5$, when the $(cp)^\sig_\tau$ was originally introduced. This part $P_5$ is unchanged up to $Fig.9$:\\
$P_5=P_6=P_7=P_9$. Roughly speaking, $\tilde{I}_0-I_1'-I_3$ can be regarded as a $\Pi_3$-series $I_0-I_1-I_3$. In this way even if the new merging points grow downwards unlimitedly, we can find a subsequence $\{\sig_{n_i}\}$ such that $st_4(\sig_{n_{i+1}})<st_4(\sig_{n_i})$. Thus any succession of collapsings terminates in a finite number of steps.
\\
\smallskip

    Once $\Pi_4-rfl$ can be analyzed, it is not so hard to treat $\Pi_n-rfl \, (n<\ome)$ and further $\Pi_\alp-rfl$ for a given transfinite $\alp<$the least $\Pi_\alp-rfl$ ordinal.
\\
    Now is the time for turning to stability from reflection.

\section{Ordinals $\sig$ having $\sig$ stable ordinals below}
\bdf {\rm For} $\alp,\bet\in Ord$ {\rm with} $\alp<\bet$, $\alp$ {\rm is} $\bet-stable$ {\rm if}
 $L_\alp\prec_{\Sig_1}L_\bet \, \Lrarw_{df} \fal \Sig_1 A \mbox{ in } L_\alp(L_\bet\models A\Lrarw L_\alp\models A)$
\edf

The reason for this turning to stability is that $\Sig^1_2$-Comprehension Axiom is interpretable in a universe $L_\bet$ such that $L_\bet$ has $\bet$-stable ordinals.

    We consider a baby case, i.e., ordinals $\sig^+$ such that $\sig$ is $\sig^+$-stable. Here recursion theoretic facts are helpful.
\\
{\bf Facts}. (cf.\cite{A-R} and \cite{M}.) For a countabl $\sig$,
\benu
\item $\sig$ is $\Pi^1_1$-reflecting $\Lrarw \: \sig$ is $\sig^+$-stable.
\item $\Pi^1_1$ on $L_\sig=$inductive on $L_\sig=\Sig_1$ on $L_{\sig^+}$.
\eenu

Let $S^1_1$ denote a theory of ordinals $\sig^+$ and $T^1_1$ a theory of ind-reflecting ordinals.

\bdf $S^1_1$ and $T^1_1$
\benu
\item {\rm The language of} $S^1_1$ {\rm is} $\calL_0\cup\{R^{\cal A}\}\cup\{\Ome,\pi\}$. {\rm The axioms of} $S^1_1$ {\rm say that the universe} $\pi^+$ {\rm of} $S^1_1$ {\rm is} $\Pi_2${\rm -reflecting and the ordinal} $\pi$ {\rm is} $\pi^+${\rm -stable: for each} $\Sig_1 \, A \; \fal u<\pi(A(u)\rarw A^\pi(u))$ {\rm or equivalently} \\
$\fal u<\pi(A(u)\rarw \exi y<\pi(y>u\& A^y(u)))$.\\
{\rm The corresponding rule runs as follows:}

\[\infer[(stbl)]{\Gam}{\Gam,\neg(t<b<\pi\wedge A^b(t)) & t<\pi\wedge A(t),\Gam}\]

\item {\rm The language of} $T^1_1${\rm is the language of} $S^1_1$ plus $\{I_<\}$, {\rm where} $I_<$ {\rm is a ternary predicate constant: Fix an} $X${\rm -positive formula} $A\equiv A(X^+,a)$ {\rm in} 
$\calL_{0}\cup\{R^{\cal A}\}\cup\{X\}$. {\rm Let} $Mp$ {\rm denote the set of multiplicative principal numbers} $a\leq\pi$ {\rm and} $a^+$ {\rm the next admissible to} $a$. {\rm Then the intended meaning of the constant} $I_<$ {\rm is given by:}
\[\fal a\in Mp\fal b<a^+[I^a_{<b}=\bigcup_{d<b}I^a_d=\bigcup_{d<b}\{c<a:A^a(I^a_{<d},c)\}]\]

{\rm That is to say, for each} $a\in Mp,a\leq\pi$ {\rm and} $b<a^+$,  $I^a_{<b}$ {\rm is the inductively generated subset of} $a=\{c:c<a\}$ {\rm by the positive formula} $A$ {\rm on the model} $<a;+,\cdot,\ldots,R^{\cal A},\ldots>$, uniformly {\rm with respect to the multiplicative principal number} $a$.\\
{\rm The axioms of} $T^1_1$ {\rm say that the universe} $\pi^+$ {\rm is} $\Pi_2${\rm -reflecting and the axiom} $(\Pi^1_1-rfl)$:
\[\fal c<\pi[c\in I^\pi_{<\pi^+}\rarw \exi \bet\in (c,\pi)\cap Mp(c\in I^\bet_{<\bet^+})].\]
{\rm where} $c\in I^a_{<a^+}\Lrarw_{df}\exi z<a^+\, A^a(I^a_{<z},c)$.
\eenu
\edf
Then it is not hard to see that $S^1_1$ is interpretable in $T^1_1$: we can extract an interpretation from Chapter 9 in Moschovakis \cite{M}.
\\
\smallskip

Before developing a proof theory for the theory $S^1_1$, we stay the theory $T^1_1$ for a while since the latter is still a theory of reflecting ordinals and an analysis for it may be attainable from $\Pi_\alp$-reflecting. We have intuitively:
  \[\mbox{Predicative Analysis }: ID_1=\Pi_\alp\mbox{-reflecting }:T^1_1\]
and since the step from Predicative Analysis to $ID_1$ requires a new dimension, an analysis for $T^1_1$ would require a new twist too.

\[\infer[(\Pi^1_1-rfl)\, J]{}
{
 \neg(\alp<b<\pi),\fal x<b^+\neg A^b(I^b_{<x},\alp) 
& 
 \infer[(\exi)]{\exi x<\pi^+ A^\pi(I^\pi_{<x},\alp)}
              {A^\pi(I^\pi_{<\xi},\alp)}
}\]
with $\alp\in I^\pi_{<\pi^+}\equiv\exi x<\pi^+ A^\pi(I^\pi_{<x},\alp)$, etc.

First consider the easy case:\\
{\bf Case1}. $\xi<\pi$: Then the above $(\Pi^1_1-rfl)\, J$ says that $\pi$ is $\Pi_\xi$-reflecting. So define $\sig=d_\pi$ such that $\xi,\alp<\sig<\pi$ and substitute $\sig$ for the variable $b$.\\
Second the general case:\\
{\bf Case2}. $\xi\geq\pi$: Pick a $\sig=d_\pi$ as above and substitute $\sig$ for $b$. We need to compute a $\xi'$ such that $\sig\leq\xi'<\sig^+$ and resolve the $(\Pi^1_1-rfl)\, J$:

\[\infer[(cut)]{}
{
 \neg A^\sig(I^\sig_{<\xi'},\alp)
&
 \infer[(cp)^\pi_\sig \, I]{A^\sig (I^\sig _{<\xi'},\alp)}
  {
    \infer[J]{A^\pi(I^\pi_{<\xi},\alp)}
            {
              \alp\not\in I^b_{<b^+} 
            &
              \alp\in I^\pi_{<\pi^+},A^\pi(I^\pi_{<\xi},\alp)
            }
   }
}
\]

The problem is that we have to be consistent with the part

\[\infer{}
{
 \neg A^\sig(I^\sig_{<\xi'},\alp)
&
 \infer{A^\sig (I^\sig _{<\xi'},\alp)}
  {A^\pi(I^\pi_{<\xi},\alp)}
}
\]

This requires a function $F :\xi\mapsto\xi'$ such that
\bdes
\item[$(F1)$] $F$ is order preserving, and in view of {\bf Case1},
\item[$(F2)$] $F$ is identity on$<\pi$, i.e., $\xi\in dom(F)\mid\pi\Rarw F(\xi)=\xi$
\item [$(F3)$]$rng(F)<\sig^+$.
\edes
Note that, here, $dom(F)$ is a {\em proper subset} of $\{\xi\in O(\pi^1_1):\xi<\pi^+\}$ with a system $O(\Pi^1_1)$ of o.d.'s for the theory $T^1_1$. We can safely set
    \[dom(F )=\{\xi\in O(\pi^1_1):K_\pi\xi<\sig\} \]
, i.e., subdiagram $\bet<\pi$ in $\xi\in dom(F)$ is $<\sig$ since $dom(F)$ is the set of o.d.'s that may occur in the upperpart of the $(cp)^\pi_\sig\, I$. Especially we have
    \[dom(F )\mid\pi=O(\pi^1_1)\mid\sig\]

   Can we take the function $F$ as a collapsing function, e.g., $d_\pi$? The answer is no. We cannot expect for $\xi,\zeta\in dom(F)$, that $\xi<\zeta \Rarw \xi\ll_\pi\zeta$ or something like an essentially less than relation. And what is worse is that the function $F$ have to preserve atomic sentences in $\calL_0$.
\bdes
\item[$(F4)$] $F$ preserves atomic sentences in $\calL_0$, i.e., diagrams of $\calL_0$ models 
     $< dom(F ) ; +,\cdot,\ldots>$ and $< rng(F ) ; +,\cdot,\ldots>$. 
\edes

To sum up $(F1)-(F4)$,
\bdes
\item[(*)] $F$ is an embedding from $\calL_0$ models $<dom(F);+,\cdot,\ldots>$ 
to $<rng(F); +,\cdot,\ldots>$ over $O(\Pi_1^1)\mid\sig$.
\edes

Now our solution for $F$ is a trite one: a {\em substitution} $[\pi:=\sig]$.
\bdes
\item [$(F5)$]$F(\xi)=\xi$ if $\xi<\pi\, (\Lrarw \xi<\sig)$
\item [$(F6)$]$F$ commutes with $+$ and the Veblen function $\vphi$, e.g., $F(\xi+\zeta)=F(\xi)+F(\zeta)$.
\item [$(F7)$]$F(\pi)=\sig$ and $F(\pi^+)=\sig^+$.
\item [$(F8)$] $F(d_{\pi^+}\bet)=d_{\sig^+}F(\bet)$.
\edes
Assume $\pi<\xi<\pi^+$ with a strongly critical $\xi$. Such a $\xi$ is of the form $d_{\pi^+}\bet$ and is introduced when a $(\Pi_2-rfl)$ for the universe $\pi^+$ is resolved. Then this $F$ meets (*), i.e., $(F1)$: Note that we have
\bdes
\item[$(F9)$] $F(K_{\pi^+}\bet)=K_{\sig^+}F(\bet)$,
\edes
and by definition $d_{\pi^+}\bet<d_{\pi^+}\gam \, \Lrarw \, 1. \, \bet<\gam \,\&\, K_{\pi^+}\bet< d_{\pi^+}\gam \mbox{ or } 2. \, d_{\pi^+}\bet\leq K_{\pi^+}\gam$ 
and similarly for $\sig^+$.
\smallskip

In this way we can resolve a $(\Pi_1^1-rfl)$ by setting $\xi'=F(\xi)$: each o.d. $\xi$ in the uppersequent of a $(cp)$ is replaced by $F(\xi)$ in the lowersequent.
\smallskip

    Next consider the theory $S^1_1$.

\[\infer[(stbl)\, J]{}
 {
  \neg(\alp<b<\pi\wedge A^b(\alp))
 &
  \infer[(\exi)]{\alp<\pi\wedge A^{\pi^+}(\alp)}
   {B(\xi,\alp)}
 }
\]
with $A^{\pi^+}(\alp)\equiv\exi x<\pi^+B(x,\alp)$.\\
As in $T_1^1$, pick a $\sig=d_\pi$ and the substitution $F=[\pi:=\sig]$. Substitute $\sig^+$ for $b$ and $\xi'=F(\xi)$ for $\xi$.

\[\infer{}
 {
  \neg B(\xi',\alp)
 &
  \infer{B(\xi',\alp)}
   {
     \neg A^{\sig^+}(\alp)
   &
     \infer[(cp)^\pi_\sig \, I]{A^{\sig^+}(\alp),B(\xi',\alp)}
                           {A^{\pi^+}(\alp),B(\xi,\alp)}
   }
 }
\]

When a universal quantifier $\fal y$ occurs in $B$, then it must be a bounded one, say, $\fal y<\xi'+\alp$ since $B$ is a bounded formula. An instance$<\xi'+\alp$ for the dual existential quantifier $\exi y<\xi'+\alp$ may come from the upperpart of $\neg B(\xi',\alp)$. Then an inspection shows that the instance$\in rng(F)$, i.e., is of the form $\zeta'=F(\zeta)$ for some $\zeta\in dom(F)$. Hence we substitute $\zeta$ for the variable $y$ in the upperpart of the $(cp)^\pi_\sig\, I$.

    In this way we can proceed and resolve consistently by (*).
\smallskip

Next we consider an ordinal $\pi$ which has many $\pi$-stable ordinals below. For example let $\pi\ome$ be an $\ome$ limit of $\pi\ome$-stable ordinals:
\[\pi\ome=\sup\{\pi n:n<\ome\} \: \& \: \fal n<\ome(\pi n \mbox{ is } \pi\ome\mbox{-stable})\]
The corresponding rule runs as follows:

\[\infer[(stbl)_n]{\Gam}{\Gam,\neg(\alp<b<\pi n\wedge A^b(\alp)) & \alp<\pi n\wedge A^{\pi\ome}(\alp),\Gam}\]

Assume $\alp<\pi n\wedge A^{\pi\ome}(\alp)$ is a conclusion of an $(\exi)$ with an auxiliary formula $B(\xi,\alp)$ with $A^{\pi\ome}(\alp)\equiv \exi x<\pi\ome B(x,\alp)$. As above we substitute 
$\xi'=F(\xi)$ for $\xi$ with $F=[\pi n:=\sig]$ for a $\sig=d_{\pi n}\bet<\pi n$ with 
$\bet=o(\alp<\pi n\wedge A^{\pi\ome}(\alp))$. 

This $F$ have to mirror the situation of o.d.'s above $\pi n$, at least occurring above the right uppersequent $\alp<\pi n\wedge A^{\pi\ome}(\alp)$. Therefore we introduce (or better postulate the existence of ordinals corresponding to) o.d.'s 
$\sig m=F(\pi m)<\pi n, \: \sig<\sig m$ for $\ome\geq m>n$. This o.d. $\sig m$ is a substitute for $\pi m$ and so have to act as if it were $\pi m$. Further when we resolve a rule $(stbl)_m$ with $m>n$, we introduce a $\tau=d_{\pi m}\gam$ with $\pi n<\tau<\pi m$ and $\tau k<\pi m$ for $\ome\geq k>m$. Thus we also have to introduce $\tau'=d_{\sig m}\gam'=F(\tau)<\sig m$ and $\tau' k=F(\tau k)<\sig m$. Then $\sig<\tau'<\tau' k<\sig m$. Let $O(2;\ome)$ denote the system of o.d.'s constructed in this way.

Here the consistency of the reduction procedure is not so problematic: these newly introduced o.d.'s are mirror images by the mirror $F$. Although $\sig m$ have to act as if it were $\pi m$, there need not be introduced a rule which says that $\sig m$ is $\sig\ome$-stable. Hence as in $S^1_1$ each instance term for an existential quantifier in $\neg B(\xi',\alp)$ is in $rng(F)$.

Rather the well foundedness of $O(2;\ome)$ is problematic: consider a series $\{\rho_i'\}$ such that $\rho_0'=\sig m$ with $n<m<\ome$, and for each $i>0$, $\rho_i'=\tau_i' (m+i)$ with $\tau_i'=d_{\rho_{i-1}'}\bet_i$ for some $\bet_i$. Then we would have a ascending sequence followed by a descending sequence:
\[\sig<\tau_1'<\cdots<\tau_k'<\tau_{k+1}'<\cdots<\rho_{k+1}'<\rho_k'<\cdots<\rho_1'<\rho_0'=\sig_m<\pi n\]
These o.d.'s came from the right upper part $\alp<\pi n\wedge A^{\pi\ome}(\alp)$ of the rule $(stbl)_n$ as mirror images by $F$. First of all preimages $\{\rho_i\}$ of these were introduced and then these are introduced as $\rho_i'=F(\rho_i)$. These preimages were created to resolve the rule $(stbl)_m$ and hence they were situated above the rule $(stbl)_n$. This means that $\rho_i<\bet=o(\alp<\pi n\wedge A^{\pi\ome}(\alp))$ and, in fact a stronger $\rho_i\ll_{\pi n^+}\bet$ holds. Therefore if we are in a situation that the o.d. $\bet$ is secured, i.e., is in a well founded part of a subrelation of $<$, then so were the descending sequence $\rho_i$. This contradicts the well foundedness.
\\
In this way we can prove that o.d.'s are well founded.
\smallskip

The whole argument works for the general case when we replace the order type $\ome$ of stable ordinals by any ordinal. Thus we get a system of o.d.'s which represent a combinatorial complexity of proof figures in a theory for ordinals $\sig$ having many $\sig$-stables. From this we also get an upper bound for the proof theoretic ordinal of a second order arithmetic for an iterated $\Sig^1_2-CA$.

\benu
\item $\Sig^1_2-CA_0$: The corresponding ordinal is a limit of ordinals $\pi n, \, n<\ome$ such that each $\pi n$ is a limit of recursively regular ordinals and has $n$ $\pi n$-stable ordinals below. A system $O(2;<\ome)$ of o.d.'s suffices.
\item $\Sig^1_2-CA+BI$: The ordinal $\pi \ome$ is a limit of $\pi\ome$-stable ordinals, i.e., {\em nonprojectible ordinal}. The set theory $KP\ome+\Sig_1 \, Separation$ is equivalent to this. $O(2;\ome)$ suffices.
\item $\Sig^1_3-DC_0$: The ordinal is a limit of ordinals $\pi a, \, a<\ome^\ome$ such that each $\pi a$ has $a$ $\pi a$-stable ordinals below. $O(2;<\ome^\ome)$ suffices.
\item $\Sig^1_3-DC$: The ordinal is a limit of ordinals $\pi a, \, a<\veps_0$ such that each $\pi a$ has $a$ $\pi a$-stable ordinals below. $O(2;<\veps_0)$ suffices.
\item $\Sig^1_3-DC+BI$: This is included in the set theory 
$KP\ome+\Pi_1 \, Collection+V=L$. 

Let $S_I$ denote a theory of ordinals $I$ such that $I$ is $\Pi_2(St^+)$-reflecting, where $St$ denotes the set of stable ordinals below $I$ and $\Pi_2(St^+)$ the set of $\Pi_2$ formulae $A$ in the language $\calL_0\cup\{St\}$ so that the predicate constant $St$ occurs only positively in the formula $A$. Then the set theory $KP\ome+\Pi_1 \, Collection+V=L$ is interpretable in $S_I$. A system $O(2;I)$ is designed for $S_I$. In $O(2;I)$ a constructor $(I,\alp)\mapsto d_I\alp\in St$ generates $I$-stable ordinals.
\eenu
Each of these systems of o.d.'s is shown to be best possible. For example we have 
\[\mid\Sig^1_2-CA+BI\mid=\mid KP\ome+\Sig_1 \, Separation\mid =O(2;\ome)\mid\Ome, \, etc.\]

\end{document}